\documentclass{amsart}
\usepackage{amssymb}
\newtheorem{thm}{Theorem}
\newtheorem{lem}{Lemma}
\newtheorem{prop}{Proposition}
\newtheorem{cor}{Corollary}

\newtheorem{rem}{Remark}
\newtheorem{df}{Definition}

\begin{document}

\bibliographystyle{plain}

\title[Milo\v s Arsenovi\'c and Romi F. Shamoyan]{On embeddings, traces and multipliers in harmonic function spaces}

\author[]{Milo\v s Arsenovi\' c$\dagger$}
\author[]{Romi F. Shamoyan}

\address{Department of mathematics, University of Belgrade, Studentski Trg 16, 11000 Belgrade, Serbia}
\email{\rm arsenovic@matf.bg.ac.rs}

\address{Bryansk University, Bryansk Russia}
\email{\rm rshamoyan@yahoo.com}

\thanks{$\dagger$ Supported by Ministry of Science, Serbia, project OI174017}

\date{}

\begin{abstract}
This paper is devoted to certain applications of classical Whitney decomposition of the upper half space
$\mathbb R^{n+1}_+$ to various problems in harmonic function spaces in the upper half space. We obtain sharp new assertions
on embeddings, distances and traces for various spaces of harmonic functions. New sharp theorems on multipliers for harmonic function spaces in the unit ball are also presented.
\end{abstract}

\maketitle

\footnotetext[1]{Mathematics Subject Classification 2010 Primary 42B15, Secondary 46E35.  Key words
and Phrases: harmonic functions, distances, traces, embedding theorems, multipliers.}

\section{Introduction, preliminaries and auxiliary results}

The role of Whitney decomposition of an open set $\Omega \subset \mathbb R^n$ in Analysis is well known, it has many
applications, for example in the theory of singular integral operators, see \cite{Gr1}, \cite{St}. The aim of this work is to present new applications of Whitney decomposition of the upper half space. In the first section we introduce notations we use in this paper and collect various auxiliary results. The second section is devoted to embedding theorems for spaces of harmonic functions in the upper half space. The third section contains results on trace problems and distance estimates and in the last one we turn from the upper half space to the unit ball $\mathbb B$ in $\mathbb R^n$ and characterize multipliers between certain spaces of harmonic functions on $\mathbb B$. Most of our results rely on Lemma \ref{LemmaA}, which provides
a Whitney type decomposition of the upper half space, and Lemma \ref{LemmaB}.

We set $\mathbb R^{n+1}_+ = \{(x, t) : x \in \mathbb R^n, t > 0 \} \subset \mathbb R^{n+1}$. For $z = (x, t) \in \mathbb R^{n+1}_+$ we set $\overline z = (x, -t)$. We denote the points in $\mathbb R^{n+1}_+$ usually by $z = (x, t)$ or $w = (y, s)$. The Lebegue measure is denoted by $dm(z) = dz = dx dt$ or $dm(w) = dw = dy ds$. We also use standard weighted measures
$dm_\lambda(z) = t^\lambda dxdt$, $\lambda \in \mathbb R$. Given two positive quantities $A$ and $B$, we write $A \asymp B$
if there are two constants $c, C > 0$ such that $cA \leq B \leq CA$.

The space of all harmonic functions in a domain $\Omega$ is denoted by $h(\Omega)$. The space of all functions $f(z_1,
\ldots, z_m)$ on $(\mathbb R^{n+1}_+)^m$ which are harmonic in each of variables $z_1, \ldots, z_m$ is denoted by
$\tilde h((\mathbb R^{n+1}_+)^m)$. Weighted harmonic Bergman spaces on $\mathbb R^{n+1}_+$ are defined, for $0<p<\infty$
and $\lambda > -1$, by
$$A^p_\lambda = A^p_\lambda(\mathbb R^{n+1}_+) = \left\{ f \in h(\mathbb R^{n+1}_+) : \| f \|_{A^p_\lambda} =
\left( \int_{\mathbb R^{n+1}_+} |f(z)|^p dm_\lambda(z) \right)^{1/p} < \infty \right\},$$
these spaces are complete metric spaces and for $1 \leq p < \infty$ they are Banach spaces.

For $\overrightarrow{s} = (s_1, \ldots, s_m)$, where $s_j > -1$, and $0<p<\infty$ we set
$$A^p_{\overrightarrow{s}} = h(\mathbb (R^{n+1}_+)^m) \cap L^p((\mathbb R^{n+1}_+)^m, dm_{s_1}(z_1) \ldots dm_{s_m}(z_m)).$$

For $p = \infty$ and $s_j > 0$, $1 \leq j \leq m$, we define $A^\infty_{\overrightarrow{s}}$ as the space of all $f \in h((\mathbb R^{n+1}_+)^m)$ such that
$$\| f \|_{A^\infty_{\overrightarrow{s}}} = \sup_{z_j \in \mathbb R^{n+1}_+} |f(z_1, \ldots, z_m)| t_1^{s_1}
\ldots t_m^{s_m} < \infty,$$
for $m=1$ and $s_1 = \alpha$ we use simpler notation $A^\infty_\alpha$. Again,
the spaces $A^p_{\overrightarrow{s}}$ are complete metric spaces, for $1 \leq p \leq \infty$ they are Banach spaces.

Finally, if $X$ is a space of functions harmonic on $(\mathbb R^{n+1}_+)^m$, then we set $\tilde hX = X \cap \tilde h ((\mathbb R^{n+1}_+)^m)$.

\begin{df}
For a function $f : (\mathbb R^{n+1}_+)^m \to \mathbb C$ we define ${\rm Tr} f : \mathbb R^{n+1}_+ \to \mathbb C$ by
${\rm Tr} f(z) = f(z,\ldots, z)$.

Let $X \subset h(\mathbb R^{n+1}_+)^m$. The trace of $X$ is ${\rm Trace}\; X = \{ {\rm Tr}\, f : f \in X \}$.
\end{df}

The problem of characterizing traces of various analytic spaces received much attention, for some results and further
references see \cite{SM1},\cite{SM2},\cite{SM3}.

We denote the Poisson kernel for $\mathbb R^{n+1}_+$ by $P(x, t)$, i.e.
$$P(x, t) = c_n \frac{t}{(|x|^2 + t^2)^{\frac{n+1}{2}}}, \qquad x \in \mathbb R^n, \quad t > 0.$$
For $l \in \mathbb N_0$ a Bergman kernel $Q_l(z, w)$, where $z = (x, t) \in \mathbb R^{n+1}_+$ and $w = (y, s) \in \mathbb R^{n+1}_+$, is defined by
$$Q_l(z, w) = \frac{(-2)^{l+1}}{l!} \frac{\partial^{l+1}}{\partial t^{l+1}} P(x-y, t+s).$$
The following theorem from \cite{DS} gives a Bergman type integral representation for functions in $A^p_\alpha$ spaces and
justifies the above terminology.
\begin{thm}\label{brthm}
Let $0<p<\infty$ and $\alpha > -1$. If $0<p\leq 1$ and $l \geq \frac{\alpha +n +1}{p} - (n+1)$ or $1\leq p < \infty$
and $l> \frac{\alpha +1}{p} - 1$, then
\begin{equation}\label{brep}
f(z) = \int_{\mathbb R^{n+1}_+} f(w) Q_l(z, w) s^l dy ds, \qquad f \in A^p_\alpha,\quad z \in \mathbb R^{n+1}_+.
\end{equation}
\end{thm}

The following elementary estimate of this kernel is contained in \cite{DS}:
\begin{equation}\label{estq}
|Q_l(z, w)| \leq C |z - \overline w|^{-(l+n+1)}, \qquad z = (x, t), \; w = (y, s) \in \mathbb R^{n+1}_+.
\end{equation}

Next we formulate two lemmata which are at the core of proofs of some of our main results. The first one provides
the above mentioned Whitney type decomposition of the upper half space. The second one is based on subharmonic
behavior of $|f|^p$ for harmonic $f$ and $0<p<\infty$.

\begin{lem}[\cite{St}]\label{LemmaA}
There exists a collection $\{ \Delta_k \}_{k=1}^\infty$ of closed cubes in $\mathbb R^{n+1}_+$ with sides parallel to coordinate axes such that

$1^o$. $\cup_{k=1}^\infty \Delta_k = \mathbb R^{n+1}_+$ and ${\rm diam} \Delta_k \asymp {\rm dist} (\Delta_k,
\partial \mathbb R^{n+1}_+)$.

$2^o$.  The interiors of the cubes $\Delta_k$ are pairwise disjoint.

$3^o$.  If $\Delta_k^\ast$ is a cube with the same center as $\Delta_k$, but enlarged 5/4 times, then the
collection $\{ \Delta_k^\ast \}_{k=1}^\infty$ forms a finitely overlapping covering of $\mathbb R^{n+1}_+$, i.e. there is
a constant $C = C_n$ such that $\sum_k \chi_{\Delta_k^\ast} \leq C$.
\end{lem}

\begin{lem}[\cite{D1}]\label{LemmaB}
Let $\Delta_k$ and $\Delta_k^\ast$ be the cubes from the previous lemma and let $(\xi_k, \eta_k)$ be the center of $\Delta_k$. Then, for $0<p<\infty$ and $\alpha > 0$, we have
\begin{equation}
\eta_k^{\alpha p - 1} \max_{\Delta_k} |f|^p \leq \frac{C}{|\Delta_k^\ast|} \int_{\Delta_k^\ast} t^{\alpha p - 1}
|f(x,t)|^p dx dt, \qquad f \in h(\mathbb R^{n+1}_+), \quad k \geq 1.
\end{equation}
\end{lem}

We will also use the following three technical estimates.

\begin{lem}[\cite{St}]\label{LemmaC}
Let $\Delta_k$ and $\Delta_k^\ast$ are as in the previous lemma, let $\zeta_k = (\xi_k, \eta_k)$ be the center of the cube $\Delta_k$. Then we have:
\begin{equation}\label{mlam}
 m_\lambda(\Delta_k) \asymp \eta_k^{n+1+\lambda} \asymp m_\lambda (\Delta_k^\ast), \qquad \lambda \in
\mathbb R,
\end{equation}
\begin{equation}\label{dis}
|\overline w - z | \asymp | \overline \zeta_k - z|, \qquad w \in \Delta_k^\ast, \quad z \in \mathbb R^{n+1}_+,
\end{equation}
\begin{equation}\label{height}
t \asymp \eta_k, \qquad (x,t) \in \Delta_k^\ast.
\end{equation}
\end{lem}

\begin{lem}[\cite{AS1}]\label{qm}
For $\delta > -1$, $\gamma > n + 1 + \delta$ and $m \in \mathbb N_0$ we have
$$\int_{\mathbb R^{n+1}_+} |Q_m(z, w)|^{\frac{\gamma}{n+m+1}} s^\delta dy ds \leq C t^{\delta - \gamma + n + 1},
\qquad z = (x, t) \in \mathbb R^{n+1}_+.$$
\end{lem}

\begin{lem}[\cite{KY}]\label{omit}
If $\alpha > -1$ and $n+\alpha<2\gamma-1$, then
\begin{equation}
\int_{\mathbb R^{n+1}_+} \frac{t^\alpha dz}{|z - \overline w|^{2\gamma}} \leq C s^{\alpha + n + 1 - 2\gamma}, \qquad
w = (y, s) \in \mathbb R^{n+1}_+.
\end{equation}
\end{lem}


\section{Embedding theorems in spaces of harmonic functions in the upper half space}

In this section we provide, using Whitney decomposition, various extensions of some embedding results from \cite{KK}.
In particular, we prove embedding theorems for $B^{p,q}_\alpha$ and $F^{p,q}_\alpha$ mixed norm classes.

Our first result relates Carleson type condition to the trace operator, as a consequence we obtain a generalization of
an embedding result from \cite{KK}.

\begin{thm}\label{tremb}
Let $0<p<\infty$, $s_1, \ldots, s_m > -1$ and let $\mu$ be a positive Borel measure on $\mathbb R^{n+1}_+$. Then the
following conditions are equivalent:

$1^o$. The measure $\mu$ satisfies a Carleson type condition:

\begin{equation} \label{cacon3}
\frac{\mu(\Delta_k)}{|\Delta_k|^{m + \frac{1}{n+1}\sum_{j=1}^m s_j}} \leq C, \qquad k \geq 1.
\end{equation}

$2^o$. ${\rm Trace}\; A^p_{\overrightarrow{s}} \hookrightarrow L^p(\mathbb R^{n+1}_+, d\mu)$, i.e. the trace
operator is bounded from $A^p_{\overrightarrow{s}}$ to $ L^p(\mathbb R^{n+1}_+, d\mu)$.

$3^o$. ${\rm Trace}\; hA^p_{\overrightarrow{s}} \hookrightarrow L^p(\mathbb R^{n+1}_+, d\mu)$, i.e. the trace
operator is bounded from $hA^p_{\overrightarrow{s}}$ to $ L^p(\mathbb R^{n+1}_+, d\mu)$.

$4^o$. If $f_j \in A^p_{s_j}$ for $1\leq j \leq m$, then
\begin{equation}\label{prod}
\int_{\mathbb R^{n+1}_+} \prod_{j=1}^m |f_j(z)|^p d\mu(z) \leq C \prod_{j=1}^m \| f_j \|_{A^p_{s_j}}^p.
\end{equation}
\end{thm}

{\it Proof.}  Let us show $1^o \Rightarrow 2^o$. Let us choose $f \in A^p_{\overrightarrow{s}}$. We use partition of $\mathbb R^{n+1}_+$ into cubes $\Delta_k$ centered at $\zeta_k = (\xi_k, \eta_k)$. Using Lemma \ref{LemmaA} we obtain
\begin{align*}
\| {\rm Tr}\; f \|_{L^p(\mu)}^p & = \int_{\mathbb R^{n+1}_+} |f(z, \ldots, z)|^p d\mu(z) = \sum_{k=1}^\infty \int_{\Delta_k} |f(z, \ldots, z)|^p d\mu(z)\\
& \leq \sum_{k=1}^\infty \mu(\Delta_k) \sup_{z \in \Delta_k} |f(z, \ldots, z)|^p
\end{align*}
Let us fix $z \in \Delta_k$. Let $B_k(z) \subset (\mathbb R^{n+1}_+)^m$ be the ball centered at $(z, \ldots, z)$ with radius equal $1/8$ of the side length of $\Delta_k$. Since $|B_k| \asymp |\Delta_k|^m \asymp \eta_k^{m(n+1)}$ and $B_k(z) \subset (\Delta_k^\ast)^m$ we obtain, using Lemma \ref{LemmaB}:
\begin{align*}
|f(z, \ldots, z)|^p & \leq \frac{C}{|B_k(z)|} \int_{B_k(z)} |f(w_1, \ldots, w_m)|^p dm(w_1)\ldots dm(w_m)\\
& \leq  \frac{C}{\eta_k^{m(n+1)}} \int_{(\Delta_k^\ast)^m} |f(w_1, \ldots, w_m)|^p dm(w_1)\ldots dm(w_m).
\end{align*}
Since the last estimate is valid for every $z \in \Delta_k$ we obtain, using (\ref{height}):
\begin{align*}
\| {\rm Tr}\; f \|_{L^p(\mu)}^p & \leq C \sum_{k=1}^\infty \frac{\mu(\Delta_k)}{\eta_k^{m(n+1)}} \int_{(\Delta_k^\ast)^m} |f(w_1, \ldots, w_m)|^p dm(w_1)\ldots dm(w_m)\\
& \leq C \sum_{k=1}^\infty \frac{\mu(\Delta_k)}{\eta_k^{m(n+1)+ \sum_{j=1}^m s_j}} \int_{(\Delta_k^\ast)^m} |f(w_1, \ldots, w_m)|^p dm_{s_1}(w_1)\ldots dm_{s_m}(w_m)
\end{align*}
Now condition (\ref{cacon3}) and finite overlapping property of $\Delta_k^\ast$ combine to give
\begin{align*}
\| {\rm Tr}\; f \|_{L^p(\mu)}^p & \leq C \sum_{k=1}^\infty \int_{(\Delta_k^\ast)^m} |f(w_1, \ldots, w_m)|^p dm_{s_1}(w_1)\ldots dm_{s_m}(w_m)\\
& \leq C \int_{(\mathbb R^{n+1}_+)^m} |f(w_1, \ldots, w_m)|^p dm_{s_1}(w_1)\ldots dm_{s_m}(w_m),
\end{align*}
which implies $2^o$.

The implication $2^o \Rightarrow 3^o$ being trivial, we turn to $3^o \Rightarrow 4^o$. Let us choose $f_j \in A^p_{s_j}$, $j=1,\ldots, m$. Using Fubini's theorem we see that $f(z_1, \ldots, z_m) = f_1(z_1) \cdots f_m(z_m)$ is in $hA^p_{\overrightarrow{s}}$ with $\| f \|_{A^p_{\overrightarrow{s}}} = \prod_{j=1}^m \| f_j \|_{A^p_{s_j}}$ and (\ref{prod}) follows immediately from $3^o$.

Finally, we prove $4^o \Rightarrow 1^o$. Let us fix a cube $\Delta_k$ and set $f_j(z) = |z - \overline \zeta_k|^{-n+1}$ for $1 \leq j \leq m$. Clearly $f_j \in A^p_{s_j}$ and, by Lemma \ref{omit}, $\| f_j \|_{A^p_{s_j}}^p \leq C \eta_k^{s_j + n + 1 - p(n-1)}$. Since $\prod_{j=1}^m f_j(z) = |z-\overline \zeta_k|^{m(1-n)}$ we have, using Lemma \ref{LemmaC},
$$\mu(\Delta_k) \eta_k^{mp(1-n)} \leq C \| \prod_{j=1}^m f_j \|_{L^p(\mu)}^p \leq C \prod_{j=1}^m \| f_j \|_{A^p_{s_j}}^p
\leq C \eta_k^{mp(1-n) + m(n+1) + \sum_{j=1}^m s_j},$$
which gives (\ref{cacon3}) and completes the proof. $\Box$


\begin{cor}\label{emTh1}
Assume $0<p<\infty$, $\alpha > -1$ and let $\mu$ be a positive Borel measure on $\mathbb R^{n+1}_+$. Then the
following conditions are equivalent.

$1^o$. The measure $\mu$ satisfies a Carleson type condition:

\begin{equation}\label{cacon1}
\frac{\mu(\Delta_k)}{|\Delta_k|^{1 + \frac{\alpha}{n+1}}} \leq C, \qquad k \geq 1.
\end{equation}

$2^o$. $A^p_\alpha$ is continuously embedded into $L^p(\mathbb R^{n+1}_+, d\mu)$.
\end{cor}

The above result was proved, in the special case $\alpha = 0$ and $1 \leq p < \infty$, in \cite{KK}.

For a measurable function $f(x, t)$ defined on $\mathbb R^{n+1}_+$ we define
$$M_p(f, t) = \| f(\cdot, t) \|_{L^p(dx)}, \qquad t > 0, \quad 0<p\leq \infty.$$

The following spaces of harmonic functions were considered, for $0<p<\infty$, $0 < q < \infty$ and $\alpha > 0$, in \cite{K1}:
$$B^{p,q}_\alpha = \left\{ f \in h(\mathbb R^{n+1}_+) : \| f \|_{B^{p,q}_\alpha} = \left( \int_0^\infty M^p_q(f, t) t^{\alpha p -1} dt \right)^{1/p} < \infty \right\},$$
in fact one can consider these spaces for $p = \infty$ or $q = \infty$, see \cite{K1}. These spaces have obvious (quasi)-norms, with respect to these (quasi)-norms they are Banach spaces or complete metric spaces.

Let us set, for $w \in \mathbb R^{n+1}_+$ and $l \geq 0$,
\begin{equation}\label{testf}
f_{w,l}(z) = \frac{\partial^l}{\partial t^l} \frac{1}{|z-\overline w|^{n-1}}.
\end{equation}
This function is harmonic in $\mathbb R^{n+1}_+$ and, for every $l \geq 0$, we have
\begin{equation}\label{polyn}
f_{w,l}(z) = \frac{1}{|z-\overline w|^{n-1+l}} P_l((t+s)|z-\overline w|^{-1}),
\end{equation}
where $P_l$ is a polynomial with integer coefficients of degree $l$. The last statement follows by induction on $l$ using identity $\partial_t |z- \overline w| = (t+s)|z-\overline w|^{-1}$. In particular we have
\begin{equation}\label{estfml}
|f_{w,l}(z)| \leq C |z - \overline w|^{-n-l+1}, \qquad z, w \in \mathbb R^{n+1}_+.
\end{equation}
Using (\ref{estfml}) one derives estimate
\begin{equation}
M_p(f_{w,l}, t) \leq C (t + s)^{\frac{n}{p} - (n-1+l)}, \qquad p(n-1+l) > n,
\end{equation}
and this estimate yields
\begin{equation}\label{bpqn}
\| f_{w,l} \|_{B^{q,p}_\alpha} \leq C s^{\frac{n}{p} - (n-1+l) + \alpha}, \qquad p(n-1+l) > n + \alpha p.
\end{equation}
Next, for $w \in \mathbb R^{n+1}_+$, let $Q_w$ be the closed cube centered at $w$ with sides parallel to the coordinate axes and with side length equal to $s$. Since the polynomial $P_l$ has finitely many zeroes it is easy to show that there are
constants $\delta > 0$ and $c > 0$ such that
\begin{equation}\label{delc}
|T_w| = c |Q_w|, \qquad T_w = \{ z \in Q_w : |P_l((t+s)|z-\overline w|^{-1})| > \delta \}.
\end{equation}

\begin{lem}\label{covarg}
Let $\mu$ be a positive Borel measure on $\mathbb R^{n+1}_+$ and let $\theta > 0$. Assume $\mu(T_w) \leq C s^\theta$ for
every $w= (y,s) \in \mathbb R^{n+1}_+$. Then we have $\mu(Q_w) \leq C s^\theta$ for every $w = (y,s) \in \mathbb R^{n+1}_+$.
\end{lem}

{\it Proof.} Let us set $T^o_w = T_w \cap Q^o_w$, where $Q^o_w$ denotes the interior of the cube $Q_w$. We note that all pairs $(Q_w, T^o_w)$ are similar to each other, either by translation by a vector parallel to the boundary $\partial \mathbb R^{n+1}_+$ or by a homothecy with center on $\partial \mathbb R^{n+1}_+$. Let us fix a cube $Q_w$. Using noted similarity it is easy to conclude that $Q_w \subset \cup \{ T^o_{w'} : w' = (y', s') \in \mathbb R^{n+1}_+, s/2 \leq s' \leq 2s \}$. Now we have a finite subcover $Q_w \subset \cup_{j=1}^N T^o_{w_j}$ where $N$ does not depend on $w$ and therefore
$$\mu(Q_w) \leq \sum_{j=1}^N \mu(T^o_{w_j}) \leq C \sum_{j=1}^N s_j^\theta \leq C s^\theta. \qquad \Box$$

The above discussion is a preparation for the following embedding theorem for $B^{p,q}_\alpha$ spaces. Its analogue for
analytic functions in the unit ball appeared in \cite{Sh3}.

\begin{thm}\label{bemb}
Let $0 < p \leq q < \infty$, $\alpha > 0$ and let $\mu$ be a Borel measure on $\mathbb R^{n+1}_+$. Then the following
conditions are equivalent:

$1^o$. The measure $\mu$ satisfies the following condition:
\begin{equation}\label{cacon2}
\mu(\Delta_k) \leq C \eta_k^{n\frac{q}{p} + \alpha q}, \qquad k \geq 1,
\end{equation}
where $(\xi_k, \eta_k)$ is the center of $\Delta_k$.

$2^o$. We have continuous embedding $B^{q,p}_\alpha \hookrightarrow L^q(\mathbb R^{n+1}_+, d\mu)$.
\end{thm}

{\it Proof.} Let us assume (\ref{cacon2}) is satisfied. We have $B^{q,p}_\alpha \hookrightarrow B^{q,q}_{\alpha + \frac{n}{p} - \frac{n}{q}} = A^q_\theta$, where $\theta = \alpha q + \frac{q}{p}n -n-1$, for the above embedding see Lemma 5 from \cite{K1}. Now we use Corollary \ref{emTh1} to obtain $A^q_\theta \hookrightarrow L^q(\mathbb R^{n+1}_+, d\mu)$.

Now assume $2^o$ is satisfied. Let us choose $l \in \mathbb N_0$ such that $p(n-1+l)>n+\alpha p$. Then, using (\ref{delc}), (\ref{bpqn}) and (\ref{polyn}), we have
\begin{equation}\label{twice}
\mu(T_w) s^{(1-n-l)q}  \leq C \| f_{w,l} \|^q_{L^q(d\mu)} \leq C \| f_{w,l} \|^q_{B^{q,p}_\alpha}
\leq C s^{q[\frac{n}{p} - (n-1+l)q] + \alpha q},
\end{equation}
which gives $\mu(T_w) \leq C s^{q \frac{n}{p} + \alpha q}$. However, Lemma \ref{covarg} shows that the
last estimate is equivalent to (\ref{cacon2}). $\Box$

In the next theorem we consider harmonic Triebel-Lizorkin spaces $F^{p,q}_\alpha$ consisting of all functions
$f \in h(\mathbb R^{n+1}_+)$ satisfying
$$ \| f \|_{F^{p,q}_\alpha}^p = \int_{\mathbb R^n} \left( \int_0^\infty |f(x,t)|^q t^{\alpha q -1} dt \right)^{p/q} dx
< \infty.$$

These spaces are complete metric spaces, for $\min (p, q) \geq 1$ these spaces are Banach spaces.

The following theorem has a counterpart for analytic functions in the unit ball in $\mathbb C^n$, see \cite{Sh3}.

\begin{thm}\label{femb}
Let $0<\tau \leq p < \infty$. Assume $\mu$ is a Borel measure on $\mathbb R^{n+1}_+$. Then the following conditions are equivalent:

$1^o$. The measure $\mu$ satisfies a Carleson type condition:
\begin{equation}\label{cacon4}
\mu(\Delta_k) \leq C \eta_k^{n + \alpha p}, \qquad k \geq 1,
\end{equation}
where $(\xi_k, \eta_k)$ is the center of $\Delta_k$.

$2^o$. The space $F^{p, \tau}_\alpha$ is continuously embedded into $L^p(\mathbb R^{n+1}_+, \mu)$.
\end{thm}

{\it Proof.} The sufficiency of condition (\ref{cacon4}) follows from the following chain of embeddings
$$F^{p, \tau}_\alpha \hookrightarrow F^{p,p}_\alpha = B^{p,p}_\alpha = A^p_\theta \hookrightarrow
L^p(\mathbb R^{n+1}, d\mu), \qquad \theta = \alpha p - 1,$$
for the first one see \cite{K1} and for the last one see Corollary \ref{emTh1}.

Next we prove necessity of the condition (\ref{cacon4}). Let us choose $l \in \mathbb N_0$ satisfying $p(n-1+l) > n + \alpha p$ and as test functions consider $f_{w,l}$ from (\ref{testf}). Our condition on $l$ ensures that $f_{w,l} \in F^{p, \tau}_\alpha$, $w = (y,s) \in \mathbb R^{n+1}_+$. Clearly the norm $\| f_{w,l} \|_{F^{p,\tau}_\alpha}$ depends on $s$ but not
on $y$, hence we can assume $w = (0, s)$, $s > 0$. A simultaneous change of variables $x = \lambda \xi$, $t = \lambda r$ where
$\xi \in \mathbb R^n$, $r > 0$, gives the following identity:
$$\| f_{(0,s),l}\|_{F^{p,\tau}_\alpha}^p = \lambda^{n-p(n-1+l-\alpha)} \| f_{(0,s/\lambda),l} \|_{F^{p,\tau}_\alpha}^p, \qquad
\lambda > 0, s>0.$$
This identity gives the following norm evaluation:
\begin{equation*}\label{normev}
\| f_{w,l} \|_{F^{p,\tau}_\alpha}^p = \| f_{(0,s),l} \|_{F^{p,\tau}_\alpha}^p = \| f_{(0,1),l} \|_{F^{p,\tau}_\alpha}^p s^{n-p(n-1+l-\alpha)}, \qquad w = (y, s) \in \mathbb R^{n+1}_+.
\end{equation*}
Therefore, we obtain the following inequalities, analogous to (\ref{twice}):
\begin{equation*}
\mu(T_w) s^{(1-n-l)p}  \leq C \| f_{w,l} \|^p_{L^p(d\mu)} \leq C \| f_{w,l} \|^p_{F^{p, \tau}_\alpha}
\leq C s^{-p(n-1+l - \alpha) +n}.
\end{equation*}
This gives $\mu(T_w) \leq C s^{n+\alpha p}$ which is, by Lemma \ref{covarg}, equivalent to (\ref{cacon4}). $\Box$

\begin{rem}
We note that in Theorems \ref{bemb} and \ref{femb} it is tempting to use $|z- \overline w|^{-n+1}$ as test functions,
they have useful property of being positive. But this works only for certain values of $\alpha$ and $p$. The problem is that this test function is not always in the required space. If it happens to be, the proof can be simplified. But in general we take derivatives and improve "size" of our test functions. However, taking derivatives in principle destroys positivity, the functions $f_{w,l}$ can have zeros and this makes proofs technically more complicated and explains the need for
Lemma \ref{covarg}.
\end{rem}

\section{Traces and distances in spaces of harmonic functions in the upper half space}

In this section, using Whitney decomposition, we present complete analogues of previously known results on traces in
analytic classes in polydisk and polyball \cite{SM1}, \cite{SM2}, \cite{SM3} and complement our previous results on distances \cite{AS2}, \cite{AS1}.

\begin{thm}
Let $s_1, \ldots, s_m > -1$ and set $\lambda = \sum_{j=1}^m s_j$. Then
\begin{equation}\label{eqtr1}
A^\infty_\lambda \subset {\rm Trace}\; \tilde h A^\infty_{\overrightarrow{s}}.
\end{equation}
Conversely, if $f \in A^\infty_{\overrightarrow{s}}$ and if ${\rm Tr}\; f$ is harmonic, then ${\rm Tr}\; f \in A^\infty_\lambda$.
\end{thm}

{\it Proof.} The second part of theorem follows from definitions. Let $g \in A^\infty_\lambda$. We choose a non negative integer $k > \lambda -1$ and define
$$f(z_1, \ldots, z_m) = \int_{\mathbb R^{n+1}_+} Q_k \left(\frac{z_1 + \cdots + z_m}{m}, w \right) g(w) s^k dw,
\qquad z_j \in \mathbb R^{n+1}_+.$$
It is immediate from Theorem \ref{brthm} that ${\rm Tr}\; f = g$. Since the kernel $Q_k$ is harmonic in
the variables $z_1, \ldots, z_m$ it is clear that $f \in \tilde h((\mathbb R^{n+1}_+)^m)$. Using classical inequality between arithmetic and geometric mean we obtain
\begin{align*}
|f(z_1, \ldots, z_m)| & \leq \| g \|_{A^\infty_\lambda} \int_{\mathbb R^{n+1}_+} \left| Q_k \left(\frac{z_1 + \cdots + z_m}{m}, w \right) \right| s^{k-\lambda} dw\\
& \leq C \| g \|_{A^\infty_\lambda} \int_0^\infty s^{k-\lambda} \int_{\mathbb R^n} \frac{dy}{|\frac{z_1 + \cdots +
z_m}{m} - \overline w|^{k+n+1}} ds\\
& \leq C \| g \|_{A^\infty_\lambda} \int_0^\infty \frac{s^{k-\lambda} ds}{(\frac{t_1 + \cdots + t_m}{m} + s)^{k+1}}\\
& \leq C \| g \|_{A^\infty_\lambda} \int_0^\infty
\frac{s^{k - \sum_{j=1}^m s_j}}{ \left( \prod_{j=1}^m (t_j + s) \right)^{(k+1)/m}}ds\\
& = C \| g \|_{A^\infty_\lambda} \int_0^\infty \prod_{j=1}^m \frac{s^{k/m - s_j}}{(t_j + s)^{\frac{k+1}{m}}} ds.
\end{align*}
Next we use H\"older's inequality for $m$ functions to obtain
\begin{equation*}
|f(z_1, \ldots, z_m)| \leq C \left( \prod_{j=1}^m \int_0^\infty \frac{s^{k-ms_j}}{(t_j + s)^{k+1}} ds \right)^{1/m}
\leq C t_1^{-s_1} \ldots t_m^{-s_m},
\end{equation*}
and the proof is finished. $\Box$

For any two $m$-tuples ($m \geq 1$) of reals $a = (a_1, \ldots, a_m)$ and $b = (b_1, \ldots, b_m)$ we define an
integral operator
\begin{equation}
(S_{a,b}f)(z_1, \ldots, z_m) = \prod_{j=1}^m t_j^{a_j} \int_{\mathbb R^{n+1}_+}
\frac{f(w) s^{-n-1+\sum_{j=1}^m b_j}}{\prod_{j=1}^m |z_j - \overline w|^{a_j + b_j}}dw, \qquad z_j \in \mathbb R^{n+1}_+.
\end{equation}
This operator is used in norm estimates of extension operators used in the proofs of trace theorems, it is well defined for $f(w) \in L^1(\mathbb R^{n+1}_+, s^{-n-1-\sum_{j=1}^m b_j})$. A unit ball analogue of this operator was used in \cite{SM1}, see also \cite{Zh}.

The following proposition is well known in the case $m=1$, its analogues for $m>1$ in polydiscs and polyballs were proved in \cite{SM1}, \cite{SM2}. The proof relies on Lemma \ref{LemmaB} and on Whitney decomposition, i.e. on Lemma \ref{LemmaA}.

\begin{prop}\label{Pr2}
Let $0<p\leq 1$, $a, b \in \mathbb R^m$ and $s_1, \ldots, s_m > -1$ satisfy $pa_j > -1-s_j$ and $pb_j>n+1+s_j$
for $j=1,\ldots, m$. Set $\lambda = (m-1)(n+1) + \sum_{j=1}^m s_j$. Then there is a constant $C>0$ such that
\begin{equation}\label{eqPr2}
\int_{\mathbb R^{n+1}_+} \cdots \int_{\mathbb R^{n+1}_+} |(S_{a,b}f)(z_1,\ldots,z_m)|^p dm_{s_1}(z_1) \ldots
dm_{s_m}(z_m) \leq C \| f \|^p_{A^p_\lambda}
\end{equation}
for every $f \in A^p_\lambda(\mathbb R^{n+1}_+)$.
\end{prop}

{\it Proof.} We use again family $\Delta_k$ of cubes from Lemma \ref{LemmaA}. We have, using (\ref{dis}),
\begin{align*}
|(S_{a,b}f)(z_1, \ldots, z_m)| & = \prod_{j=1}^m t_j^{a_j} \sum_{k=1}^\infty \int_{\Delta_k}
\frac{f(w) s^{-n-1+\sum_{j=1}^m b_j}}{\prod_{j=1}^m |z_j - \overline w|^{a_j + b_j}}dw\\
& \leq C \sum_{k=1}^\infty \prod_{j=1}^m \frac{t_j^{a_j}}{|z_j - \overline\zeta_k|^{a_j+b_j}} \int_{\Delta_k} |f(w)|
s^{-n-1+\sum_{j=1}^m b_j} dm(w)\\
& \leq C \sum_{k=1}^\infty \prod_{j=1}^m \frac{t_j^{a_j}\eta_k^{\sum_{j=1}^m b_j}}{|z_j - \overline\zeta_k|^{a_j+b_j}}
\sup_{\Delta_k} |f|.
\end{align*}
Since $0<p\leq 1$ this gives
\begin{equation}\label{toint}
|(S_{a,b}f)(z_1, \ldots, z_m)|^p \leq C \sum_{k=1}^\infty \prod_{j=1}^m
\frac{t_j^{pa_j}\eta_k^{p\sum_{j=1}^m b_j}}{|z_j - \overline\zeta_k|^{pa_j+pb_j}} \sup_{\Delta_k} |f|^p.
\end{equation}
We integrate this inequality with respect to $dm_{s_1}(z_1) \ldots dm_{s_m}(z_m)$ and obtain, using Lemma \ref{omit}
and Lemma \ref{LemmaB}:
\begin{align*}
M & = \int_{\mathbb R^{n+1}_+} \cdots \int_{\mathbb R^{n+1}_+} |(S_{a,b}f)(z_1,\ldots,z_m)|^p
dm_{s_1}(z_1) \ldots dm_{s_m}(z_m)\\
& \leq C \sum_{k=1}^\infty \eta_k^{m(n+1)+\sum_{j=1}^m s_j} \sup_{\Delta_k} |f|^p  \leq C \sum_{k=1}^\infty \int_{\Delta_k^\ast} |f(w)|^p dm_\lambda(w).
\end{align*}
This is sufficient due to finite overlapping property of the family $\Delta_k^\ast$. $\Box$

\begin{lem}\label{Lemma4}
Let $0<p<\infty$ and $s_1, \ldots, s_m > -1$. Set $\lambda = (m-1)(n+1) + \sum_{j=1}^m s_j$. Then there is a constant
$C>0$ such that for all $f \in h((\mathbb R^{n+1}_+)^m)$ we have
\begin{align}\label{eqL4}
\int_{\mathbb R^{n+1}_+} |f(z,\ldots, z)|^p & dm_\lambda (z) \notag \\
& \leq C \int_{\mathbb R^{n+1}_+} \cdots \int_{\mathbb R^{n+1}_+} |f(z_1, \ldots, z_m)|^p dm_{s_1}(z_1) \ldots dm_{s_m}(z_m).
\end{align}
\end{lem}

{\it Proof.} This is an immediate consequence of Theorem \ref{tremb}, indeed it is easy to check that measure
$dm_\lambda$ satisfies condition (\ref{cacon3}). $\Box$

\begin{rem}
We note that the trace results from \cite{SM1} for harmonic spaces in the unit ball are true under one additional condition: namely that the trace of a considered harmonic function on a product of unit balls is harmonic on the unit ball.
\end{rem}

The theorem below was announced, without proof, in \cite{SM1} for $0<p<\infty$. Analogous results in the case of the unit ball were proved in \cite{SM1}. Let us note that in the case $n=1$ a precise subclass for which we have characterization of
traces was found in \cite{SM3}, this is precisely the space of all pluriharmonic functions in the unit poly disk.

\begin{thm}\label{tr1T}
Let $0<p \leq 1$, $s_1, \ldots, s_m > -1$ and set $\lambda = (m-1)(n+1) + \sum_{j=1}^m s_j$. Then
\begin{equation}\label{eqtr1}
A^p_\lambda \subset {\rm Trace}\; \tilde h A^p_{\overrightarrow{s}} \subset {\rm Trace}\; A^p_{\overrightarrow{s}} \subset
L^p(\mathbb R^{n+1}_+, dm_\lambda).
\end{equation}
In particular, if $f \in A^p_{\overrightarrow{s}}$ and if ${\rm Tr}\; f$ is harmonic,
then ${\rm Tr}\; f \in A^p_\lambda$.
\end{thm}

{\it Proof.} The second inclusion in (\ref{eqtr1}) is trivial and the last one follows from Lemma \ref{Lemma4}, which also establishes the last part of the theorem. It remains to prove $A^p_\lambda \subset {\rm Trace}\; \tilde h A^p_{\overrightarrow{s}}$. Let
$g \in A^p_\lambda$ and set
$$f(z_1, \ldots, z_m) = \int_{\mathbb R^{n+1}_+} Q_l\left(\frac{z_1 + \cdots + z_m}{m}, w \right) g(w) s^l dw,
\qquad z_j \in \mathbb R^{n+1}_+,$$
where $l > \lambda - 1$ is an integer. As in the previous theorem, clearly $f \in \tilde h((\mathbb R^{n+1}_+)^m)$. In order to simplify notation we set $Z = z_1 + \cdots + z_m$. Now we have, using (\ref{height}), (\ref{estq}) and (\ref{dis}),
\begin{align*}
|f(z_1, \ldots, z_m)| & \leq \sum_{k=1}^\infty \int_{\Delta_k} |Q_l(Z/m, w)| \; |g(w)| s^l dw\\
& \leq C \sum_{k=1}^\infty \eta_k^l \int_{\Delta_k} |Q_l(Z/m, w)| dw \sup_{\Delta_k} |g|\\
& \leq C \sum_{k=1}^\infty \frac{\eta_k^{n+l+1}}{|Z/m - \overline \zeta_k|^{n+1+l}} \sup_{\Delta_k} |g|.
\end{align*}
Since $0<p\leq 1$ this estimate and Lemma \ref{LemmaB} give
\begin{align*}
|f(z_1, \ldots, z_m)|^p & \leq C \sum_{k=1}^\infty \frac{\eta_k^{p(n+l+1)}}{|Z/m - \overline \zeta_k|^{p(n+1+l)}} \sup_{\Delta_k} |g|^p\\
& \leq C \sum_{k=1}^\infty \frac{\eta_k^{p(n+l+1)-\lambda-n-1}}{|Z/m - \overline \zeta_k|^{p(n+1+l)}} \| g \|^p_{L^p(
\Delta_k^\ast, dm_\lambda)}.
\end{align*}
Now we set $z_j = (x_j, t_j)$, $1\leq j \leq m$ and obtain
\begin{align*}
I_k(t_1, \ldots, t_m) & = \int_{\mathbb R^n} \cdots \int_{\mathbb R^n}
\frac{dx_1 \ldots dx_m}{\left| \frac{(x_1 + \cdots + x_m, t_1 + \cdots
t_m)}{m} - \overline \zeta_k \right|^{p(n+l+1)}}\\
& \leq C \left( \frac{t_1 + \cdots + t_m}{m} + \eta_k \right)^{-p(n+l+1) + mn}.
\end{align*}
Inequality between arithmetic and geometric mean and the above estimate imply
\begin{align*}
M & = \int_{\mathbb R^{n+1}_+} \cdots \int_{\mathbb R^{n+1}_+} |f(z_1, \ldots, z_m)|^p dm_{s_1}(z_1) \ldots
dm_{s_m}(z_m)\\
& \leq C \sum_{k=1}^\infty \eta_k^{p(n+l+1) - \lambda-n-1} \int_0^\infty \cdots \int_0^\infty
\frac{t_1^{s_1} \ldots t_m^{s_m}dt_1 \ldots dt_m}{ \left( \frac{t_1 + \cdots t_m}{m} + \eta_k \right)^{p(n+l+1) - mn}}
\| g \|^p_{L^p(\Delta_k^\ast, dm_\lambda)}\\
& \leq C \sum_{k=1}^\infty \eta_k^{p(n+l+1) - \lambda-n-1} \prod_{j=1}^m \int_0^\infty
\frac{t_j^{s_j} dt_j}{(t_j + \eta_k)^{\frac{p(n+l+1)-mn}{m}}} \| g \|^p_{L^p(\Delta_k^\ast, dm_\lambda)}\\
& \leq C \sum_{k=1}^\infty \| g \|^p_{L^p(\Delta_k^\ast, dm_\lambda)} \leq C \| g \|_{A^p_\lambda}
\end{align*}
due to finite overlapping property of the family $\Delta_k^\ast$. $\Box$

\begin{rem}
As already noted in \cite{SM2}, the above results on traces can be extended to certain mixed norm spaces of harmonic functions.
\end{rem}

We note that $A^p_\alpha$ is continuously embedded into $A^\infty_{\frac{\alpha + n + 1}{p}}$, see \cite{K1}. Hence it is
a natural problem to look for estimates of ${\rm dist}_{A^\infty_{\frac{\alpha + n + 1}{p}}} (f, A^p_\alpha)$ for $f \in
A^\infty_{\frac{\alpha+n+1}{p}}$. The case $1<p<\infty$ was treated in \cite{AS1}, the theorem below covers the remaining
case $0<p\leq 1$.

For $\epsilon > 0$, $\lambda > 0$ and $f \in h(\mathbb R^{n+1}_+)$ we set
$$V_{\epsilon, \lambda} (f)= V_{\epsilon, \lambda} = \left\{ (x, t) \in \mathbb R^{n+1}_+ : |f(x,t)|t^\lambda \geq \epsilon \right\}.$$

\begin{thm}\label{displ1}
Let $0 < p \leq 1$, $\alpha > -1$, $\lambda = \frac{\alpha + n + 1}{p}$, $m \in \mathbb N_0$ and $m > \max(\lambda -1, \frac{\alpha}{p})$. Set, for $f \in A^\infty_{\frac{\alpha + n +1}{p}}$,
\begin{equation}
d_1(f) = {\rm dist}_{A^\infty_{\frac{\alpha + n + 1}{p}}} (f, A^p_\alpha),
\end{equation}
\begin{equation}\label{d2f}
d_2(f) = \inf \left\{ \epsilon > 0 : \int_{\mathbb R^{n+1}_+} \left( \int_{V_{\epsilon, \lambda}} |Q_m(z, w)|
s^{m - \lambda}
dy ds \right)^p t^\alpha dx dt < \infty \right\}.
\end{equation}

Then $d_1(f) \asymp d_2(f)$.
\end{thm}

{\it Proof.} We begin with inequality $d_1(f) \geq d_2(f)$. Assume $d_1(f) < d_2(f)$. Then there are $\epsilon > \epsilon_1 > 0$ and $f_1 \in A^p_\alpha (\mathbb R^{n+1}_+)$ such that $\| f - f_1 \|_{A^\infty_\lambda} \leq \epsilon_1$ and
\begin{equation}\label{dpl1i}
\int_{\mathbb R^{n+1}_+} \left( \int_{V_{\epsilon, \lambda}} |Q_m(z, w)| s^{m - \lambda} dy ds \right)^p t^\alpha dx dt = + \infty.
\end{equation}
Since  $\| f - f_1 \|_{A^\infty_\lambda} \leq \epsilon_1$ we have
\begin{equation}\label{dpl1e}
(\epsilon - \epsilon_1) \chi_{V_{\epsilon, \lambda}}(w) s^{-\lambda} \leq |f_1(w)|, \qquad w = (y, s) \in \mathbb
R^{n+1}_+.
\end{equation}
Combining (\ref{dpl1i}) and (\ref{dpl1e}) we obtain
\begin{align*}
+ \infty & = \int_{\mathbb R^{n+1}_+} \left( \int_{\mathbb R^{n+1}_+} \chi_{V_{\epsilon, \lambda}}(w)
|Q_m(z, w)| s^{m - \lambda} dy ds \right)^p t^\alpha dx dt \\
& \leq C \int_{\mathbb R^{n+1}_+} \left( \int_{\mathbb R^{n+1}_+} |f_1(w) Q_m(z, w)| s^m dy ds \right)^p t^\alpha
dx dt = M.
\end{align*}
Our goal is to obtain a contradiction by showing that $M$ is finite. We use cubes $\Delta_k$ and $\Delta_k^\ast$
from Lemma \ref{LemmaA} and Lemma \ref{LemmaB}, let $\zeta_k = (\xi_k, \eta_k)$ be the corresponding centers. Using Lemma
\ref{LemmaC}, Lemma \ref{LemmaA} and assumption $0<p\leq 1$ we obtain

\begin{align*}
I(z) & = \left( \int_{\mathbb R^{n+1}_+}|f_1(w) Q_m(z, w)| s^m dy ds \right)^p  = \left( \sum_{k=1}^\infty \int_{\Delta_k} |f_1(w) Q_m(z,w)| s^m dy ds \right)^p\\
& \leq C \left( \sum_{k=1}^\infty \eta_k^m |\Delta_k| \max_{w \in \Delta_k} |f_1(w)| \max_{w\in\Delta_k} |Q_m(z,w)|
\right)^p\\
& \leq C \sum_{k=1}^\infty \eta_k^{mp} |\Delta_k|^p \max_{w\in\Delta_k} |f_1(w)|^p \max_{w\in\Delta_k} |Q_m(z,w)|^p.
\end{align*}
Therefore, using Lemma \ref{LemmaB}, Lemma \ref{LemmaC} and finite overlapping property of the family $\Delta_k^\ast$ we obtain
\begin{align*}
I(z) & \leq C \sum_{k=1}^\infty \eta_k^{mp} |\Delta_k|^p |z - \overline \zeta_k|^{-p(m+n+1)} \eta_k^{1-\alpha p} \frac{1}{|\Delta_k^\ast|} \int_{\Delta_k^\ast} s^{\alpha p -1} |f_1(w)|^p dy ds\\
& \leq C \sum_{k=1}^\infty \eta_k^{mp} |\Delta_k|^{p-1} \int_{\Delta_k^\ast} \frac{|f_1(w)|^p dy ds}{|z-\overline w|^{p(m
+n+1)}}\\
& \leq C \sum_{k=1}^\infty \int_{\Delta_k^\ast}
\frac{|f_1(w)|^p s^{mp + (n+1)(p-1)} dy ds}{|z-\overline w|^{p(m+n+1)}}\\
& \leq C \int_{\mathbb R^{n+1}_+} \frac{|f_1(w)|^p s^{mp + (n+1)(p-1)} dy ds}{|z-\overline w|^{p(m+n+1)}}.
\end{align*}
This estimate, Lemma \ref{omit} with $2\gamma = p(m+n+1)$ and Fubini's theorem yield
\begin{align*}
M & \leq C \int_{\mathbb R^{n+1}_+} \left(\int_{\mathbb R^{n+1}_+} \frac{|f_1(w)|^p s^{mp + (n+1)(p-1)} dw}{|z-\overline w|^{p(m+n+1)}} \right)
t^\alpha dz \\
& = C \int_{\mathbb R^{n+1}_+} |f_1(w)|^p s^{mp + (n+1)(p-1)} \int_{\mathbb R^{n+1}_+} \frac{t^\alpha dz}{|z -
\overline w|^{p(m+n+1)}} dw\\
& \leq C \int_{\mathbb R^{n+1}_+} |f_1(w)|^p s^\alpha dw < \infty,
\end{align*}
arriving at a contradiction. Therefore we proved $d_1(f) \geq d_2(f)$.

Next we prove $d_1(f) \leq C d_2(f)$. We choose $\epsilon > 0$ such that the integral appearing in (\ref{d2f}) is finite. Using Theorem \ref{brthm} we obtain
\begin{equation*}
f(z)  = \int_{\mathbb R^{n+1}_+ \setminus V_{\epsilon, \lambda}} f(w) Q_m(z,w) s^m dw + \int_{V_{\epsilon, \lambda}}
f(w) Q_m(z,w)s^m dw  = f_1(z) + f_2(z).
\end{equation*}
Since the kernel $Q_m(z,w)$ is harmonic in the second variable both $f_1$ and $f_2$ are harmonic in $\mathbb R^{n+1}_+$.
Next, using Lemma \ref{qm} and definition of the set $V_{\epsilon, \lambda}$ we obtain
\begin{equation}
|f_1(z)| \leq \epsilon \int_{\mathbb R^{n+1}_+} |Q_m(z,w)| s^{m-\lambda} dw \leq C\epsilon t^{-\lambda}, \qquad
z = (x, t) \in \mathbb R^{n+1}_+,
\end{equation}
which gives $\| f_1 \|_{A^\infty_\lambda} \leq C\epsilon$. To complete the proof it suffices to show that $f_2 \in
A^p_\alpha$. Since $|f(w)| \leq \| f \|_{A^\infty_\lambda} s^{-\lambda}$ we have
\begin{equation*}
\| f_2 \|_{A^p_\alpha}  \leq \| f \|_{A^\infty_\lambda} \int_{\mathbb R^{n+1}} \left( \int_{V_{\epsilon, \lambda}}
|Q_m(z, w)| s^{m-\lambda} dw \right)^p t^\alpha dz \leq C  \| f \|_{A^\infty_\lambda},
\end{equation*}
and the proof is complete. $\Box$

\section{Multipliers on spaces of harmonic functions in the unit ball}

In this section we prove new sharp theorems on multipliers in spaces of harmonic functions in the unit ball
$\mathbb B = \{ x \in \mathbb R^n : |x| < 1 \}$. This topic is fairly new, it was initiated in \cite{SA} and pursued further in \cite{AS1}, \cite{AS3}. This section presents new results in this direction. Even the case of the unit disc was not studied extensively, for some results we refer to \cite{Pa2}. Of course, the topic of multipliers between analytic function spaces is
a vast subject.

We denote spherical harmonics of order $k$ by $Y^{(k)}_j$, $1 \leq j \leq d_k$, see \cite{StW} for details on spherical
harmonics. Let us recall, for reader's convenience, needed definitions.

\begin{df}
For a double indexed sequence of complex numbers
$$c = \{ c_k^j : k \geq 0, 1 \leq j \leq d_k \}$$
and a harmonic function $f(rx') = \sum_{k=0}^\infty r^k \sum_{j=1}^{d_k} b_k^j(f) Y^{(k)}_j(x')$ we define
$$(c \ast f) (rx') = \sum_{k=0}^\infty \sum_{j=1}^{d_k} r^k c_k^j b_k^j(f) Y^{(k)}_j(x'), \qquad rx' \in \mathbb B,$$
if the series converges in $\mathbb B$. Similarly we define convolution of $f, g \in h(\mathbb B)$ by
$$(f \ast g)(rx') = \sum_{k=0}^\infty \sum_{j=1}^{d_k} r^k b_k^j(f)b_k^j(g) Y_j^{(k)}(x'), \qquad rx'\in \mathbb B,$$
it is easily seen that $f \ast g$ is defined and harmonic in $\mathbb B$.
\end{df}

\begin{df}
Let $X$ and $Y$ be subspaces of $h(\mathbb B)$. We say that a double indexed sequence $c$ is a multiplier from $X$ to $Y$ if $c \ast f \in Y$ for every $f \in X$. The vector space of all multipliers from $X$ to $Y$ is denoted by $M_H(X, Y)$.
\end{df}

We are looking for sufficient and/or necessary condition for a double indexed sequence $c$ to be in $M_H(X, Y)$, for certain spaces $X$ and $Y$ of harmonic functions. We associate to such a sequence $c$ a harmonic function
\begin{equation}\label{gc}
g_c(x) = g(x) = \sum_{k\geq 0} r^k \sum_{j=1}^{d_k} c_k^j Y^{(k)}_j(x'), \qquad x = rx' \in \mathbb B,
\end{equation}

The conditions we are looking for are expressed in terms of fractional derivatives of $g = g_c$.

\begin{df}
For $t \in \mathbb R \setminus \mathbb Z_{-}$ and a harmonic function $f(x) = \sum_{k=0}^\infty r^k b_k(f)\cdot Y^k(x')$ on $\mathbb B$ we define a fractional derivative of order $t$ of $f$ by the following formula:
$$(\Lambda_t f)(x) = \sum_{k=0}^\infty r^k \frac{\Gamma(k+n/2 + t)}{\Gamma(k+n/2)\Gamma(t)}b_k(f)\cdot Y^k(x'),
\qquad x = rx' \in \mathbb B.$$
\end{df}

It is easily seen that $\Lambda_t$ maps $h(\mathbb B)$ into $h(\mathbb B)$.

We record the following formula
\begin{equation}\label{simple}
(c \ast f)(r^2x') = \int_{\mathbb S} (g_c \ast P_{y'})(rx') f(ry') dy' = \int_{\mathbb S} (g_c \ast P_{x'})(ry') f(ry') dy'
\end{equation}
where $P_{y'}(x) = P(x, y') = c_n \frac{1-|x|^2}{|x-y'|^n}$ is the Poisson kernel for the unit ball $\mathbb B$, see \cite{AS1} for details and further references.





Next we recall definitions of some harmonic function spaces. For $\alpha > 0$ we set $A^\infty_\alpha = \{ f \in
h(\mathbb B) : \| f \|_{A^\infty_\alpha} = \sup_{x \in \mathbb B} (1-|x|^2)^\alpha |f(x)| < \infty \}$ and
$$A^p_\alpha = \left\{ f \in h(\mathbb B) : \| f \|_{A^p_\alpha} = \left( \int_{\mathbb B} |f(x)|^p(1-|x|^2)^\alpha dx
\right)^{1/p} < \infty \right\}, \qquad 0<p<\infty.$$
For $0<p<\infty$, $0 \leq r < 1$ and $f \in h(\mathbb B)$ we set
$$M_p(f, r) = \left( \int_{\mathbb S} |f(rx')|^p dx' \right)^{1/p},$$
with the usual modification to cover the case $p = \infty$. The harmonic Hardy spaces are defined by
$H^s = \{ f \in h(\mathbb B): \sup_{0\leq r < 1} M_s(f, r) < \infty \}$, $0<s\leq \infty$. For $0 < p \leq \infty$,
$0 < q \leq \infty$, $\alpha > 0$ and $f \in h(\mathbb B)$ we consider mixed (quasi)-norms $\| f \|_{p,q,\alpha}$ defined by
\begin{equation}\label{qpnorm}
\| f \|_{p,q,\alpha} = \left( \int_0^1 M_q(f, r)^p (1-r^2)^{\alpha p - 1} r^{n-1}dr \right)^{1/p}, \qquad 0<p<\infty,
\end{equation}
and, for $p = \infty$, $\| f \|_{\infty, q, \alpha} = \sup_{0 \leq r < 1} (1-r^2)^\alpha M_q(f, r)$.
The corresponding mixed norm spaces are
$$B^{p,q}_\alpha(\mathbb B) = B^{p,q}_\alpha = \{ f \in h(\mathbb B) : \| f \|_{p,q,\alpha} < \infty \}.$$
For details on these spaces we refer to \cite{DS}, Chapter 7. In particular these spaces are complete metric spaces and for $\min(p,q) \geq 1$ they are Banach spaces.

Our next result uses duality arguments. For $M > -1$, $0<p<\infty$ and $\alpha > -1$ we denote by $D_{-M}A^p_\alpha$ the space of all functions $f \in h(\mathbb B)$ such that
$$\| f \|_{p,M;\alpha} = \left( \int_{\mathbb B} |\Lambda_{M+1}f(x)|^p (1-|x|^2)^\alpha dx \right)^{1/p} < \infty.$$
It is immediate that each $f \in D_{-M}A^1_{M-\beta}$, $M > \beta - 1$, generates a continuous linear functional on
$A^\infty_\beta$ by the following formula:
$$L_f(g) = \int_{\mathbb B} \Lambda_{M+1}f(x) g(x) (1-|x|^2)^M dx, \qquad g \in A^\infty_\beta.$$
Therefore, the above pairing gives an embedding $D_{-M}A^1_{M-\beta} \hookrightarrow (A^\infty_\beta)^\ast$, $M > \beta -1$, which is used in the proof of the next theorem.

\begin{thm}\label{harber}
Let $1<s<\infty$, $\beta > 0$ and let $s'$ be the exponent conjugate to $s$.
Then $c \in M_H(H^s, A^\infty_\beta)$ if and only if the function $g = g_c$ satisfies the following condition
\begin{equation}\label{eqmul1}
N_{s'}(g) = \sup_{0\leq\rho < 1} \sup_{y' \in \mathbb S} (1-\rho)^\beta
\left( \int_{\mathbb S} |(g \ast P_{x'})(\rho y')|^{s'} dx' \right)^{1/s'} < \infty.
\end{equation}
\end{thm}

{\it Proof.} Assume $c \in M_H(H^s, A^\infty_\beta)$, i.e $M_c : H^s \rightarrow A^\infty_\beta$ and choose $M>\beta-1$. Then, using the above embedding, the adjoint operator $M_c^\ast$ maps $D_{-M}A^1_{M-\beta}$ into $H^{s'}$. Moreover, it is easy to verify that $M_c^\ast$ acts as a multiplier operator from $D_{-M}A^1_{M-\beta}$ to $H^{s'}$ generated by the same double indexed sequence $c$. Next, using definition of the space $D_{-M}A^1_{M-\beta}$ we see that the double indexed sequence
$$ c' = \left\{  \frac{\Gamma(k+n/2)\Gamma(M+1)}{\Gamma(k+n/2 + M+1)}c_j^k, k\geq 0, 1 \leq j \leq d_k \right\}$$
acts as a multiplier from $A^1_{M-\beta}$ to $H^{s'}$. Set $g' = g_{c'}$, then we have $g = \Lambda_{M+1} g'$. Since $A^1_{M-\beta} = B^{1,1}_{M-\beta+1}$ we have $c' \in M_H(B^{1,1}_{M-\beta+1}, H^{s'})$ and Theorem 1 from \cite{AS3} gives, with $m > M-\beta$:
$$\sup_{0\leq \rho < 1} \sup_{y' \in \mathbb S} (1- \rho)^{-(M-\beta+1) + m + 1} \left( \int_{\mathbb S}
|\Lambda_{m + 1}(g' \ast P_{x'})(\rho y')|^{s'} dx' \right)^{1/s'} < \infty.$$
But this implies
\begin{align*}
& \sup_{0\leq \rho < 1} \sup_{y' \in \mathbb S} (1- \rho)^{\beta + m + 1} \left( \int_{\mathbb S}
|\Lambda_{m + 1}(g \ast P_{x'})(\rho y')|^{s'} dx' \right)^{1/s'} \\
= & C \sup_{0\leq \rho < 1} \sup_{y' \in \mathbb S} (1- \rho)^{\beta + m +1} \left( \int_{\mathbb S}
|\Lambda_{m + M + 2}(g' \ast P_{x'})(\rho y')|^{s'} dx' \right)^{1/s'}\\
\leq & C \sup_{0\leq \rho < 1} \sup_{y' \in \mathbb S} (1- \rho)^{\beta + m -M} \left( \int_{\mathbb S}
|\Lambda_{m + 1}(g' \ast P_{x'})(\rho y')|^{s'} dx' \right)^{1/s'} < \infty.
\end{align*}
Since $\beta > 0$, the last estimate implies (\ref{eqmul1}).

Now we assume $c$ satisfies (\ref{eqmul1}) and choose $f \in H^s$. Set $h = M_c f$. Then, using (\ref{simple}), we have
\begin{equation*}
|h(r^2 x') | \leq \int_{\mathbb S} |(g \ast P_{x'})(ry')f(ry')| dy' \leq N_{s'}(g) (1-r)^{- \beta} M_s(f,r),
\end{equation*}
which means that $(1-r)^\beta M_\infty(h, r) \leq C \| f \|_{H^s}$. $\Box$

For the case $0<s\leq 1$ we refer reader to \cite{AS3}.

The next theorem deals with $B^{p,q}_\alpha$ spaces on the unit ball. For another results on multipliers from $B^{p,q}_\alpha$ spaces we refer to \cite{AS3}.

\begin{thm}
Let $1<p<\infty$, $1\leq q < \infty$, $\alpha, \beta > 0$ and let $q'$ be the exponent conjugate to $q$. Assume $m > \max(\alpha - \beta - 1, \beta -1)$. Then $c \in M_H(B^{p,q}_\alpha, A^\infty_\beta)$ if and only if the function $g = g_c$ satisfies the following condition
\begin{equation}\label{eqmul2}
M_{q'}(g) = \sup_{0\leq\rho < 1} \sup_{y' \in \mathbb S} (1-\rho)^{m + 1 + \beta -\alpha}
\left( \int_{\mathbb S} |\Lambda_{m+1}(g \ast P_{x'})(\rho y')|^{q'} dx' \right)^{1/q'} < \infty.
\end{equation}
\end{thm}

{\it Proof.} For necessity of this condition we use a duality result $(B^{p,q}_\alpha)^\ast \cong B^{p',q'}_\alpha$ from
\cite{Za1}, see \cite{AS3} for another application of this duality to multiplier problems. Since the arguments are
analogous to those in the proof of Theorem \ref{harber} we omit details.

Now we prove sufficiency of the condition (\ref{eqmul2}). Assume $f \in B^{\infty, q}_\alpha$ and assume $c$ satisfies
(\ref{eqmul2}). Set $h = M_c f$. Using (\ref{simple}) we obtain
\begin{align*}
| \Lambda_{m+1} h(r^2 x') | & \leq \int_{\mathbb S} | \Lambda_{m+1}(g \ast P_{x'})(ry')f(ry')| dy'\\
& \leq M_{q'}(g) (1-r)^{\alpha - \beta -m-1} M_q(f,r),
\end{align*}
which gives
$$(1-r)^{m+1+\beta} M_\infty(\Lambda_{m+1} h, r) \leq C (1-r)^\alpha M_q(f,r) \leq C \| f \|_{B^{\infty, q}_\alpha}.$$
Hence $(1-r)^\beta M_\infty(h, r) \leq C\| f \|_{B^{\infty, q}_\alpha}$, see \cite{DS} Chapter 7. Therefore $M_c : B^{\infty, q}_\alpha \rightarrow A^\infty_\beta$. Since we have embedding $B^{p,q}_\alpha \hookrightarrow B^{\infty,q}_\alpha$ for all $0<p\leq \infty$, see \cite{DS}, the proof is completed. $\Box$

Finally we state two results on multipliers in spaces where definitions of norms involve derivatives, proofs will appear
elsewhere. Let $\nabla f = (\partial f/\partial x_1, \ldots, \partial f/\partial x_n)$ denote the gradient of a smooth
function $f$ on $\Omega \subset \mathbb R^n$. We set
\begin{equation*}
DA^p_\alpha = \{ f \in h(\mathbb B) : \|f\|_{DA^p_\alpha} = |f(0)| + \| \nabla f \|_{A^p_\alpha} < \infty \}, \qquad
\alpha > 0, 0<p<\infty,
\end{equation*}
\begin{equation*}
DB^{p,q}_\alpha = \{ f \in h(\mathbb B) : \| f \|_{DB^{p,q}_\alpha} = |f(0)| + \| \nabla f \|_{B^{p,q}_\alpha} < \infty \},
\qquad \alpha > 0, 0<p,q<\infty.
\end{equation*}

\begin{thm}\label{muldb}
Let $1<s<\infty$, $\alpha, \beta > 0$, $0<p\leq 1$ and let $s'$ be the exponent conjugate to $s$.
Then $c \in M_H(DB^{p,1}_\alpha, H^s_\beta)$ if and only if the function $g = g_c$ satisfies the following condition
\begin{equation}\label{eqmul3}
L_{s'}(g) = \sup_{0\leq\rho < 1} \sup_{y' \in \mathbb S} (1-\rho)^{m +2 + \beta -\alpha}
\left( \int_{\mathbb S} |\Lambda_{m+1}(g \ast P_{x'})(\rho y')|^{s'} dx' \right)^{1/s'} < \infty.
\end{equation}
\end{thm}

Since $DA^p_\alpha = DB^{p,p}_{\frac{\alpha + 1}{p}}$, see \cite{DS}, taking $p = 1$ we obtain the following corollary.

\begin{cor}
Let $1<s<\infty$, $\alpha, \beta > 0$ and let $s'$ be the exponent conjugate to $s$. Then $c \in M_H(DA^1_\alpha, H^s_\beta)$
if and only if the function $g = g_c$ satisfies the following condition
\begin{equation}\label{eqmul4}
K_{s'}(g) = \sup_{0\leq\rho < 1} \sup_{y' \in \mathbb S} (1-\rho)^{m +1 + \beta -\alpha}
\left( \int_{\mathbb S} |\Lambda_{m+1}(g \ast P_{x'})(\rho y')|^{s'} dx' \right)^{1/s'} < \infty.
\end{equation}
\end{cor}

\end{document}